\documentclass[11pt]{amsart}
\usepackage[margin=3cm]{geometry}
\title{Mixing time bounds for edge flipping on regular graphs}

\author{Yunus Emre Demirci}
\address{Queen's University, Department of Mathematics and Statistics, Kingston, Ontario, Canada}
 \email{21yed@queensu.ca }
  \author{\"{U}m\.{i}t I\c{s}lak}\address{Bo\u{g}azi\c{c}i University, Department of Mathematics, Istanbul, Turkey} 
 \email{umit.islak1@boun.edu.tr}  
\author{Alperen  \"{O}zdemir}
\address{Department of Mathematics, Georgia Institute of Technology} 
\email{aozdemir6@gatech.edu} 

\usepackage{amssymb, amsmath, amsthm, enumerate, mathtools}
\usepackage[colorlinks=true, urlcolor=black,linkcolor=blue, citecolor=blue]{hyperref}

\newtheorem{theorem}{Theorem}[section]
\newtheorem{lemma}{Lemma}[section]
\newtheorem{corollary}{Corollary}[section]
\newtheorem{definition}{Definition}[section]
\newtheorem{example}{Example}[section]

\newcommand{\bbox}{\hfill $\Box$}
\newcommand{\pf}{\noindent {\it Proof:} }

 \newcommand{\su}{\textnormal{supp }}

\newcommand{\tn}[1]{\textnormal{#1}}


\begin{document}

\begin{abstract}
The edge flipping is a non-reversible Markov chain on a given connected graph, which is defined by Chung and Graham in \cite{CG12}. In the same paper, its eigenvalues and stationary distributions for some classes of graphs are identified. We further study its spectral properties to show a lower bound for the rate of convergence in the case of regular graphs. Moreover, we show that a cutoff occurs at $\frac{1}{4} n \log n$ for the edge flipping on the complete graph by a coupling argument.   
 \end{abstract}

\maketitle

\section{Introduction}\label{sec:intro}

 The edge flipping is a random process on graphs where an edge of a graph is randomly chosen at each step and both of its endpoints are colored either blue with probability $p$ or red with probability $q=1-p$. We are interested in the long-term vertex color configurations of the graph and a derived statistic, the frequency of blue vertices. In most cases, the initial configuration of vertex colors is taken monochromatic (all-blue or all-red). Stationary distributions for paths and cycles \cite{CG12}, and for complete graphs \cite{B15} are studied and various asymptotic results on stationary color configurations are obtained for those cases. Our emphasis will be on the mixing time of the chain.
 
The set of color configurations of a graph with $n$ vertices is in fact the $n$-dimensional hypercube, $\{0,1\}^{n}.$ So it is expected that the edge flipping has similarities with other Markov chains on the space. Take, for example, the Ehrenfest urn model. In this model, there are two urns containing $n$ balls in total. At every stage a ball is randomly chosen and is moved to the opposite urn. In our case, the urns can be associated with two different colors and the dynamics is picking two balls at random and placing them together in one of the two urns. The crucial difference is that the edge flipping is non-reversible unlike the former model. In order to see that, consider the two configurations, one with all vertices blue and the other with all but one vertex blue. Then the chain moves from the latter to the former with positive probability, but the other direction has zero probability. 

The edge flipping can also be considered a random walk on the set of chambers of Boolean arrangements, which are particular cases of hyperplane arrangements. The eigenvalues of hyperplane chamber walks were first obtained in \cite{BHR99}. Markov chains such as Tsetlin library and various riffle shuffle models are shown to be examples of these walks in the same paper. The transition matrix of the random walk is shown to be diagonalizable in \cite{BD98} using algebraic topology, and a non-compelling condition for the uniqueness of the stationary distribution is provided. In a broader context, Brown in \cite{B00} gives a ring theoretical proof for diagonalizability for the random walks on left-regular bands, which include hyperplane arrangements. Finally, a combinatorial derivation of the spectrum of the random walk is shown in \cite{AD10}. The eigenvectors of hyperplane arrangements are studied in \cite{JP, Sal}.

We use the spectral analysis of the Markov chain, which is developed in the aforementioned papers, and combine it with various probabilistic techniques to address the rate of convergence to its stationary distribution with respect to total variation distance. See \cite{N19} for an example of separation distance bounds on hyperplane random walks. The spectral techniques are limited due to the fact that the edge flipping is non-reversible. In particular, a frequently used $l^2$ bound on the total variation distance by the sum of squares of eigenvalues is not available (See Lemma 12.16 in \cite{L}). An upper bound rather involving the sum of the eigenvalues is obtained by combining coupling techniques with spectral analysis in \cite{BHR99} and \cite{BD98}.

We are ready to present our results on the mixing time of the edge flipping. The mixing time of a Markov chain is the time required for the distance between the law of the chain to the stationary distribution to be small. The formal definitions for the distance notion and the mixing time are given in the beginning of Section \ref{sec:edgeflip}. Our first result is for the most general case. In Theorem \ref{thm:gen}, we show an upper and a lower bound on the mixing time for the edge flipping on connected graphs with $n$ vertices, which differ by a logarithmic factor of $n$. Then we specialize to regular graphs and obtain bounds of order $n \log n$ with different constants in Theorem \ref{epthm}. This is known as \textit{pre-cutoff}, which will not be defined here but we refer to Section 18.1 of \cite{L} for its definition. A finer result is known as \textit{cutoff},  which refers to a sharp transition in the distance of the Markov chain to its stationary distribution. We refer the reader to Section \ref{sec:complete} for its definition. 

A complete characterization of cutoff is still an open problem. There are conditions that guarantee its existence in broad classes of Markov chains. For example, \textit{the product condition} implies the existence cutoff for some classes of reversible chains \cite{BHP17}, yet there are counter-examples too, see Section 18, Notes in \cite{L}. Recently, an entropic criterion is proposed in \cite{S21} which requires \textit{the symmetry of probabilities}, a less restrictive property than reversibility. However, that condition does not apply to the edge flipping either for the aforementioned reason of it not being reversible. So, we are not able to tell immediately whether the edge flipping shows a cutoff or not. The main result of the paper is to find matching constants for the upper and the lower bound on the mixing time in the case of complete graphs to establish the following.

\begin{theorem}\label{compthm}
Let $K_n$ be the complete graph with $n$ vertices. The edge flipping on $K_n$ exhibits a cutoff at time $\frac{1}{4} n \log n$ with a window of order $n$.
\end{theorem}

In order to have an intuition for this exact rate, let us reconsider the random walks on hypercubes. In particular, let us take the lazy random walk on $\{0,1\}^n.$ In comparison to the lazy random walk on the hypercube, the number of vertices that change color at each step is $1$ on average instead of $1/2.$ So we expect to have a twice as fast chain as the former, which happens to be the case, see Theorem 18.3 in \cite{L}.


The outline of the paper is as follows. The next section provides background on hyperplane arrangements and the Markov chains defined on them. In Section 
\ref{sec:edgeflip} we formally present our results. Then we discuss how to bound the rate of convergence to the stationary distribution and supply the proofs of our theorems.

\section{Hyperplane random walks}

Let $\mathcal{A}=\{H_i\}_{i=1}^n$ be a collection of hyperplanes in a real vector space $V.$ The collection is called $central$ if $\bigcap H_i \neq \emptyset.$ Each $H_i$ divides $V$ into two open half-spaces. Let us denote them by $H_i^{+}$ and $H_i^{-},$ and let $H_i^{0}$ stand for the hyperplane itself. A \textit{face} $F$ is a non-empty subset of $V$ determined by hyperplanes as
\begin{equation*}
F= \bigcap_{i=1}^n H_i^{\sigma_i(F)}
\end{equation*}
where $\sigma_i(F) \in \{ +, -, 0 \}.$ So we can represent any face by a sign sequence $\{\sigma_i(F)\}_{i=1}^n.$ The faces with $\sigma_i(F)\neq 0$ for all $i$ are called \textit{chambers}. The set of faces and the set of chambers are denoted by $\mathcal{F}$ and $\mathcal{C}$ respectively. Then we define a left-multiplication on $\mathcal{F}$ by

\begin{equation} \label{leftmult}
\begin{aligned}
   \sigma_i(F F') = 
    \begin{cases}
      \sigma_i(F)  \text{  if  } \sigma_i(F) \neq 0,\\
      \sigma_i(F')  \text{  if  } \sigma_i(F)=0.
    \end{cases}
\end{aligned}
\end{equation}
The operation defined above generates a random walk on $\mathcal{C}$, which has transition probabilities
\begin{equation} \label{trans3}
P(C,C')= \sum_{\substack{F \in \mathcal{F} \\ FC=C'}}  w(F)
\end{equation}
where $C, C' \in \mathcal{C}$ and $w$ is any probability distribution over faces.

To state the results on the eigenvalues of the matrix associated with the transition probabilities, we first consider a partial order on $\mathcal{A}$ which is defined via \eqref{leftmult} as 
\begin{equation}\label{partord1}
F \leq F' \Leftrightarrow FF'=F'.
\end{equation}
This implies that $F \leq F'$ if and only if either $\sigma_i(F)= \sigma_i(F')$ or $\sigma_i(F)=0$ for each $i = 1,\ldots,n$. We can view the left-multiplication by $F$ as a projection from $\mathcal{A}$ onto 
$$\mathcal{A}_{\geq F} := \{F' \in \mathcal{A} : F \leq F'\}.$$ 

Next, we define another partial order by 
\begin{equation} \label{partord2}
F \preceq F' \Leftrightarrow F' F=F'.
\end{equation}
In other words, $F \preceq F'$ if and only if $\sigma_i(F')=0$ implies $\sigma_i(F)=0.$ The equivalence classes of \eqref{partord2} give a semilattice $L,$ and $L$ is a lattice if and only if $\mathcal{A}$ is central \cite{St07}.  The equivalence class of $F$ with respect to \eqref{partord2} is called the \textit{support} of $F$, and is denoted by  $\su F.$  In lattice terminology, we have 
\begin{equation*}
\su F F' = \su F \vee \su F'.
\end{equation*}
The elements of $L$ are called \textit{flats}. Furthermore, if $\su F = \su F',$ then the projection space $\mathcal{A}_{\geq F}$ is isomorphic to $\mathcal{A}_{\geq F'}.$ So the following is well-defined;
\begin{equation}
 c_X=|\mathcal{A}_{\geq F}|
 \end{equation} 
 for any $F \in \mathcal{F}$ with $X=\su F.$

Observe that the support of the chambers is maximal in $L,$ and the set of chambers is closed under left-multiplication. Therefore, the Markov chain is well-defined over the set of chambers. Its spectral profile is as follows.
\begin{theorem}\label{thm:CG} \emph{(\cite{BD98}, Theorem 1)}
Let $\{ w(F)\}_{F \in \mathcal{F}}$ be a probability distribution over the faces of a central hyperplane $\mathcal{A}$. The eigenvalues of the random walk given by transition probabilities \eqref{trans3} are indexed by the lattice $L.$ For each flat $X\in L,$ the associated eigenvalue is
\begin{equation*}
\lambda_X = \sum_{\su F \subseteq X} w(F) 
\end{equation*}
with multiplicity $m_X$ satisfying
\begin{equation*}
\sum_{X \preceq Y} m_Y = c_X.
\end{equation*}
\end{theorem}
We note that the multiplicities have the explicit expression
\begin{equation*}
m_X= \sum_{X \preceq Y} \mu(X,Y) c_Y
\end{equation*}
by M\"{o}bius inversion where $\mu$ is the M\"{o}bius function of the lattice $L.$ For a Boolean arrangement, the partial order \eqref{partord2} is the set inclusion, and its M\"{o}bius function is simply $\mu(X,Y)=(-1)^{|Y|-|X|}$ by the Inclusion-Exclusion Principle. In this case, the eigenvalues are indexed by the Boolean lattice.

Finally, we look at the edge flipping case. Let $G$ be a connected graph with $n$ vertices and $m$ edges. Denote the vertex set of $G$ by $V(G)$ and the edge set by $E(G).$ Each vertex $i$ is associated with a hyperplane $H_i$ in $\mathbb{R}^n$ and the arrangement are central. All color configurations where the vertex $i$ is blue reside in the half-space $H_i^+$ and the color configurations where the vertex $i$ is red fall into $H_i^-.$

The faces that correspond to color configurations on subsets of vertices are the chambers of the hyperplane where the sign of $\sigma_k(C)$ gives the color of the vertex $k$ in the subset. In agreement with half-spaces defined above, if $\sigma_k(F)=+ (-),$ the vertex of the subset labeled by $k$ is blue (red). The faces that generate the random walk are colored edges. A blue edge connecting the vertices $i$ and $j$ is given by the face which has the sign sequence
\begin{equation*}
   \sigma_k(F_b) = 
    \begin{cases}
      +, \text{  if  } k \in \{i,j\},\\
     \, 0, \text{  otherwise.  } 
    \end{cases}
 \end{equation*} 
In the same way, a red edge connecting $i$ and $j$ is given by
 \begin{equation*}
   \sigma_k(F_r) = 
    \begin{cases}
      -, \text{  if  } k \in \{i,j\},\\
      \, 0, \text{  otherwise.  } 
    \end{cases}
 \end{equation*} 
The probabilities assigned to these faces are uniform for the same color, which are given by $w(F_b)=\frac{p}{m}$ and $w(F_r)=\frac{q}{m}.$ 

\section{Edge flipping}\label{sec:edgeflip}

We are ready to study the stationary distribution and the rate of convergence to it for the edge flipping. We first define the distance notion to be used. The \textit{total variation distance} between $\mu$ and $\nu$ on the state space $\Omega$ is defined as
\begin{equation}\label{tv}
\| \mu - \nu \|_{TV} = \frac{1}{2} \sum_{x \in \Omega} |\mu(x) - \nu(x)| = \max_{S \subseteq \Omega}|\mu(S) - \nu(S)|.
\end{equation}
In the context of Markov chains, we will use the following notation. For a Markov chain defined on $\Omega,$ let the probability assigned to  $x \in \Omega$  be $P^{t}_z(x)$ after running the chain for $t$ steps with the initial state $z$, and let it be $\pi(x)$ at the stationary distribution. A shorthand notation will be used for the total variation distance between $P^t$ and $\pi$:
\begin{equation*}
d_{TV}(t)= \max_{z \in \Omega}\|P^t_z - \pi\|_{TV}.
\end{equation*}
In case that it does not make difference from which state we initiated the Markov chain, we will denote its law after $t$ steps by $P^t.$ In fact, it can be shown that the maximum can only be achieved by point masses on color configurations. For the rate of convergence, we define the mixing time
\begin{equation*}
t_{\text{mix}}(\epsilon)=\min \{t:  d_{TV}(t) \leq \epsilon\}.
 \end{equation*}

\subsection{Distance to stationary distribution}

The following theorem gives an upper bound to the distance to stationary distributions in terms of the eigenvalues given in the previous section.

\begin{theorem}\label{thm:BHR} \emph{(\cite{BHR99}, Theorem 5.1)}
Consider a central hyperplane arrangement with set of flats $L.$ Let $L^*$ be the set of all flats in $L$ except the unique maximal one, and $\mu$ be the M\"{o}bius function of the lattice $L.$ We have
\begin{equation*}
 d_{TV}(t) \leq - \sum_{X \in L^*} \mu(X,V(G)) \lambda_X^t
\end{equation*}
where $\lambda_X,$ the eigenvalue associated with flat $X,$ is as defined in Theorem \ref{thm:CG}. 
\end{theorem}

As mentioned in the introduction, this upper bound is not obtained by eigenvalue techniques, symmetric matrices etc. This is due to a coupling argument

Let us give two examples where we identify the eigenvalues of the trasition matrix and apply the theorem above.

\begin{example}\emph{(The complete graph, $K_n$)}
$K_n$ is the graph with $n$ vertices where every vertex is connected to every other vertex. Let $X_k$ be a flat in $L$ consisting of $k$ vertices. We have $\mu(X_k, V(G)) = (-1)^{n-k}.$
By Theorem \ref{thm:CG}, the eigenvalue corresponding to flat $X_k$ is $\frac{\binom{k}{2}}{\binom{n}{2}}$ and the number of flats of size $k$ is $\binom{n}{k}.$  Hence,
\begin{flalign*}
d_{TV}(t) &\leq -\sum_{X \in L^*} \mu(X,V(G)) \lambda_X^t \\
 &= -\sum_{k=2}^{n-1}\sum_{X_k} \mu(X_k,V(G)) \left(\frac{\binom{k}{2}}{\binom{n}{2}}\right)^t \\
 &= -\sum_{k=2}^{n-1} (-1)^{n-k} \binom{n}{k} \left(\frac{\binom{k}{2}}{\binom{n}{2}}\right)^t \\
 & \leq  n \left(1-\frac{2}{n} \right)^t +\sum_{i=2}^{n-2} (-1)^{i-1} \binom{n}{i}  \left(1-\frac{i}{n} \right)^{2t} + o(n)\\
 & = \mathcal{O}\left( n \left(1-\frac{2}{n} \right)^t\right)
\end{flalign*}
provided that $t \geq \frac{1}{2}n \log n.$ 
\end{example}

\begin{example}\emph{(The complete bipartite graph, $K_{m,n}$)}
$K_{m,n}$ is a graph whose vertex set is partitioned into two sets of sizes  $m$ and $n$ where every vertex of one set is connected to every vertex of the other set. Take a flat $X_{k,l} \in L$ where $k$ and $l$ denote the numbers of vertices from the set of size $m$ and the set of size $n$ respectively. The reader can verify using Theorem \ref{thm:CG} that $\lambda_{X_{k,l}}=\frac{kl}{mn}$ with multiplicity $\binom{m}{k} \binom{n}{l}.$ One can also show that $\mu(X_{k,l}, V(G)) = (-1)^{m+n-k-l}$. So by Theorem \ref{thm:BHR}, 
\begin{flalign*}
d_{TV}(t) &\leq -\sum_{X \in L^*} \mu(X,V(G)) \lambda_X^t \\
 &=  n \left(\frac{(n-1)m}{nm}\right)^t+ m \left(\frac{n (m-1)}{nm}\right)^t -\sum_{k=1}^{m-1}\sum_{l=1}^{n-1} (-1)^{m+n-k-l}\binom{m}{k} \binom{n}{l} \left(\frac{kl}{mn} \right)^t  \\
 &=  n \left(1- \frac{1}{n}  \right)^t + m \left(1- \frac{1}{m}  \right)^t + \sum_{j=1}^{m-1} (-1)^j \binom{m}{j} \left(1- \frac{j}{m}  \right)^t \, \sum_{i=1}^{n-1} (-1)^{i} \binom{n}{i} \left(1- \frac{i}{n}  \right)^t.  
\end{flalign*}
Suppose $n\geq m.$ Then, for $t \geq n \log n,$ we have
\[d_{TV}(t) =\mathcal{O}\left( n \left(1-\frac{1}{n} \right)^t\right).\]

\end{example}

A computationally tractable bound is obtained in \cite{BD98} by a coarser version of the coupling argument in the proof of Theorem \ref{thm:BHR} in \cite{BHR99}. It is as follows.
\begin{equation}\label{BDbound}
d_{TV}(t) \leq - \sum_{X \in L^*} \mu(X,V(G)) \lambda_X^t \leq \sum_{ \{ M: M \textnormal{is co-maximal in } L \}} \lambda_{M}^t,
\end{equation}
where $M$ is a \textit{co-maximal} flat in $L$ if $M \prec X$ implies $X$ is maximal in $L.$ In the edge flipping, the co-maximal flats are simply obtained by removing a vertex from the flat of all vertices since any flat is a subset of vertices.

The precision of the alternating bound compared to \eqref{BDbound} is discussed in \cite{C14} considering various examples of random hyperplane arrangements. In the examples that we study, the latter bound is as good as the former. 

Now we can provide an upper bound for the rate of convergence of the edge flipping in general connected graphs.  

\begin{theorem} \label{thm:gen}
Consider the edge flipping on a connected graph $G$ with $n$ vertices and $m$ edges, and let $\delta$ be the degree of the vertex with the minimum degree in $G.$ For $c>0,$ we have
\begin{equation*}
c\frac{m}{\delta} \leq t_{\text{\tn{mix}}}(qe^{-2c}) \leq \frac{m}{\delta}( \log n + 2c - \log q).
\end{equation*}
\end{theorem}

\pf  Observe that the co-maximal flats obtained by removing any vertex with the minimum number of edges have the smallest eigenvalue among all co-maximal flats. Let $M^*$ denote this flat. Then the eigenvalue indexed by $M^*$ is
\begin{equation*}
\lambda_{M^*} = \sum_{\su F \subseteq M^{*}} w(F) = \frac{m- \delta}{m}=1-\frac{\delta}{m}
\end{equation*} 
by Theorem \ref{thm:CG}. By \eqref{BDbound}, we have
$$
d_{TV}(t)  \leq \sum_{ \{ M: M \textnormal{is co-maximal in } L \}} \lambda_{M}^t \leq  n \, \lambda_{M^{*}}^t \leq  n \left( 1- \frac{\delta}{m}\right)^t \leq e^{-tm/ \delta} \leq qe^{-2c}
$$
if $t \geq \frac{m}{\delta}(\log n + 2c -\log q).$

For the lower bound, let us label a vertex with the minimal degree by $x$ and take its color blue in the initial configuration without loss of generality. Let $R$ be the set of all configurations where $x$ is colored red. We have $\pi(R)=q$ because whenever an edge with an endpoint $x$ is chosen, it is colored blue with probability $p$ and red with probability $q.$ Now let us look at $P^t(R).$ The probability that the edge connecting $x$ to the rest of the graph is not chosen before step $t$ is
\[\left(1-\frac{\delta}{m} \right)^t \geq e^{-t\left(\frac{m}{\delta} + \frac{m^2}{\delta^2}\right)}.\]
If we take $t\leq cm/ \delta,$ then we have $P^t(R) \leq q(1-e^{-2c}).$ Finally using the second inequality in \eqref{tv}, we have
\[\|P^t-\pi\|_{TV} \geq |P^t(R)-\pi(R)| \geq qe^{-2c}. \]
\bbox

We note that this bound is not optimal in either direction which is to be shown in the following section.

\subsection{Rate of convergence for regular graphs}\label{sec:edgereg}

A graph $G$ is called \textit{regular} if every vertex of $G$ has the same degree. It is called \textit{k-regular} if the degree of the vertices is $k.$ We show the following bounds on the mixing time of edge flipping for regular graphs, which are independent of the degree of the vertices.

\begin{theorem}\label{epthm}
Let $G$ be a connected, regular graph with $n$ vertices. Then for $c>0,$ the mixing time of edge flipping satisfies 
\begin{equation*}
\frac{1}{4}  n \log n - \mathcal{O}(n) \leq t_{\text{mix}}(e^{-c}) \leq \frac{1}{2}  n \log n + \mathcal{O}(n).
\end{equation*} 
\end{theorem}
The smallest degree regular graph on $n$ vertices is the cycle $C_n$ and the largest degree regular graph on $n$ vertices is $K_n$. We remark that the stationary distributions of the edge flipping for these graphs are studied in \cite{CG12} and \cite{B15}. 

 The  derivation of these bounds shows that as long as the number of vertices of the same degree is of order $n$, the principal terms in the bounds above remain the same.

\noindent \textbf{Proof of the upper bound in Theorem \ref{epthm}}.  Suppose degree of vertices of $G$ is equal to $k$. Let $M$ denote a co-maximal flat. Then the corresponding eigenvalue  is 
\begin{equation*}
\lambda_M = \sum_{\su F \subseteq M} w(F) = \frac{\frac{kn}{2}-k}{\frac{kn}{2}}=1-\frac{2}{n}
\end{equation*}
by Theorem \ref{thm:CG}. The bound \eqref{BDbound} gives
$$
\|P^t_C - \pi\|_{TV} \leq n \left(1- \frac{2}{n} \right)^t 
 \leq e^{-c}
$$
if we take $t \geq n\log n/2 +  cn/2.$ \bbox
 
For the proof of the lower bound in Theorem \ref{epthm}, we use the Wilson's method. Wilson \cite{Wilson} showed that an eigenvector is a good candidate for a test statistic as the variance of the associated eigenfunction can be estimated inductively from the transition probabilities of the Markov chain. Assuming the second-order estimate for the eigenfunction, a lower bound can be obtained as follows.

\begin{lemma} \emph{(\cite{L}, Theorem 13.5)}\label{Wilson}
Let $X_t$ be an irreducible, aperiodic, time-homogenous Markov chain with state space $\Omega.$ Let $\Phi$ be an eigenfunction associated with eigenvalue $\lambda > \frac{1}{2}.$ If for all $x \in \Omega,$
$$ \mathbf{E} \left( \left(\Phi(X_1)- \Phi(x)\right)^2 \, | \,  X_0=x \right) \leq R $$
for some $R > 0,$ then
\begin{equation*}
t_{\text{mix}}(\epsilon) \geq \frac{1}{2 \log (1 / \lambda)} \left( \log \left( \frac{(1- \lambda) \Phi(x)^2}{2R} \right) + \log \left( \frac{1- \epsilon}{\epsilon}\right)\right)
\end{equation*}
for any $x \in \Omega.$
\end{lemma}

The maximum eigenvalue is $1$ with right eigenfunction $\phi_{0}(C)=1$ for all $C \in \mathcal{C}.$ The eigenvectors of eigenvalues corresponding to co-maximal flats in hyperplane arrangements are identified by Pike as follows.

\begin{theorem}\emph{(\cite{JP}, Theorem 3.1.1)} \label{pike}
For each $i \in \{1, \dots, n \},$ the Markov chain defined on the chambers by the transition probabilities \eqref{trans3} has  the right eigenfunction
\begin{equation*} 
\begin{aligned}
   &   \phi_i(C)=
   \begin{cases}
     -\sum_{\substack{F \in \mathcal{F} \\ \sigma_i(F)=+}}  w(F),&   \sigma_i(C)=-\\
    \quad  \sum_{\substack{F \in \mathcal{F} \\ \sigma_i(F)=-}}  w(F),&   \sigma_i(C)=+
    \end{cases}
\end{aligned}
\end{equation*}
corresponding to eigenvalue $\lambda_i=\sum_{\sigma_i(F)=0} w(F)$. Moreover,  $\phi_0, \phi_1, \cdots, \phi_n$ are linearly independent where $\phi_0$ is the eigenvector for the largest eigenvalue $1.$
\end{theorem}

In our problem, the eigenvectors have a simpler expression as 
\begin{equation} \label{pikevector}
\begin{aligned}
   &   \phi_i(C)=
   \begin{cases}
     - p,&   \sigma_i(C)=-\\
    \quad  q,&   \sigma_i(C)=+
    \end{cases}
\end{aligned}
\end{equation}
for $i=1, \dots, n.$ Each $i$ represents a vertex and we can read off the color of $i$ in a chamber by the eigenfunction $\phi_i.$ These eigenfunctions are not quite distinguishing by their own. In fact, Wilson's lemma applied to any $\phi_i$ gives a lower bound of order $n$ only. Nevertheless, if an eigenvalue corresponding to a co-maximal flat has algebraic multiplicity larger than one, linear combinations of eigenvectors can be used in Wilson's lemma. In the case of regular graphs, where each vertex has the same degree, we have the second largest eigenvalue
\begin{equation*}
\lambda=\sum_{\sigma_1(F)=0} w(F)= \frac{kn/2-k}{kn/2}=1-\frac{2}{n} 
\end{equation*}
associated with co-maximal flat obtained by removing the vertex indexed by $1.$ Of course, it has multiplicity $n$ considering all other co-maximal flats and equality of degrees. 
Therefore, from Theorem \ref{pike}, the eigenspace for the second largest eigenvalue is spanned by eigenvectors \eqref{pikevector} for $i=1, \dots, n.$  If we consider the linear combination  $\Phi = \sum_{i=1}^n \phi_i,$ it is easy to see that $\Phi(C) = \Phi(C')$ if and only if they have the same number of blue vertices. If the color configuration $C$ has $k$ blue vertices, then $\Phi(C)= k -pn.$ So, $\Phi$ is a statistic that counts the number of blue vertices up to an additive term which guarantees $\mathbf{E}(\Phi(\pi))=0$ where the expectation is taken with respect to the stationary distribution $\pi,$ i.e.,
\begin{equation*}
\mathbf{E}(\Phi(X_t | X_0=C))=\lambda^{t} \Phi(C) \rightarrow 0.
\end{equation*}
 We can use this statistic to obtain the lower bound by Wilson's method. Since at each step at most two vertices can change color, we have
\begin{equation}\label{rbound}
\mathbf{E}[(\Phi(X_1)- \Phi(C))^2 | X_0=C)] \leq 4,
 \end{equation} 
so that the condition in Lemma \ref{Wilson} is satisfied with $R=4$.
 
 \vspace{0.1in}
 
\noindent \textbf{Proof of the lower bound in Theorem \ref{epthm}}.
For the lower bound, since the second eigenvalue has multiplicity $n$ given that the graph is regular, the sum of eigenvectors is an eigenvector of the second largest eigenvalues above. Take the initial state to be the blue monochromatic graph, for which $\Phi(C)= n-pn=qn$. Then using Lemma \eqref{Wilson} 
$$
 t_{\text{mix}}(\epsilon)  \geq \frac{n}{4} \left( \log \left(\frac{q^2n^2}{4n}\right) + \log \left( \frac{1- \epsilon}{\epsilon}\right)\right) 
  \geq  \frac{1}{4}n \log n -  Cn
$$
for some constant $C$ depending only on $p$ and $\epsilon$. We take $\epsilon=e^{-c}$ to conclude the proof. \bbox

 We also note that, thanks to the anonymous referee who pointed out this, the construction of eigenfunctions of edge flipping is analogous to the construction of eigenfunctions for the simple random walk on the hypercube. This also justifies the lower bound on the leading order of $\frac{1}{4} n \log n$. 
 
\subsection{Exact rate of convergence for complete graphs}\label{sec:complete}

 In Section \ref{sec:edgereg}, we established a lower bound $\frac{1}{4}n \log n$ for the mixing time of edge flipping in regular graphs with $n$ vertices. However, the upper bound obtained from eigenvalues differs by a factor of two. So we can only give an interval for the exact rate of convergence. Yet, some earlier studies on the similar processes suggest that the lower bound captures the correct rate of convergence, see \cite{N17} for instance. We can verify that in the case of complete graphs. First, we formally define what we mean by the ``exact rate."

 \begin{definition}\label{cutoffdef}
 A sequence of Markov chains shows a cutoff at the \textit{mixing time} $\{t_n\}$ with a \textit{window} of size $\{w_n\}$ if
\begin{itemize}
\item[\textnormal{(i)}] $\lim\limits_{n\rightarrow \infty} \frac{w_n}{t_n} = 0,$
\item[\textnormal{(ii)}] $\lim\limits_{\gamma \rightarrow - \infty} \liminf\limits_{n\rightarrow \infty } d_{TV}(t_n + \gamma w_n) = 1,$
\item[\textnormal{(iii)}] $\lim\limits_{\gamma \rightarrow  \infty} \limsup\limits_{n\rightarrow \infty } d_{TV}(t_n + \gamma w_n) = 0.$
\end{itemize}
\end{definition}
We will use a coupling argument to show that the edge flipping in complete graphs shows a cutoff at time $\frac{1}{4}n \log n$ with a window of size $n.$  The coupling time of the two copies of a process is defined as
\begin{equation*}
\tau:= \min\{t  :  X_t =Y_t\}.
\end{equation*}
Writing $d(\cdot)$ for $d_{T V}(\cdot)$, the coupling time is related to   the mixing time as follows:
\begin{equation}\label{coupletime}
d(t) \leq \max_{x,y \in \Omega} \mathbf{P} (\tau > t \, | \, X_0=x,Y_0=y).
\end{equation}
 See, for example, Theorem 5.4 in \cite{L}.  
\subsubsection{Proof of Theorem \ref{compthm} }

Consider the linear combination of eigenvectors denoted by $\Phi$ in Section \ref{sec:edgereg}. We already showed that it counts the number of blue vertices. Now if we view the edge flipping on $K_n$ as a random walk on the hypercube $(\mathbb{Z} / 2 \mathbb{Z})^n $ where at each step two coordinates are replaced either by 1's by probability $p$ or by 0's with probability $q,$ then $\Phi$ is just the \textit{Hamming weight}, the sum of the 0-1 coordinates. See Section 2.3. of \cite{L} for the definition of hypercube random walk. An $n \log n$ upper bound is obtained for lazy random walk on the hypercube by a coupling and a strong stationary time argument, which can be found in Section 5.3 and Section 6.4 of \cite{L}, respectively. This bound is later improved by monotone coupling in Section 18.2 of the same book following the observation that the problem can be reduced to convergence of Hamming weights, which can be considered as a lazy Ehrenfest urn problem. The same reduction argument is valid in the case of edge flipping in complete graphs, which is as follows. We claim that if the initial configuration is monochromatic, then
\begin{equation} \label{egalite}
\| P^{t} - \pi \|_{TV} = \| \Phi(X_t) - \pi_{\Phi} \|_{TV} 
\end{equation}
where $X_t$ is the chamber at step $t$ and $\pi_{\Phi}$ is the stationary distribution for the number of blue vertices. To argue for this, let us denote the set of chambers with $k$ blue vertices by $\mathcal{C}_k=\{C \in \mathcal{C} : \Phi(C)=k \}.$ Starting from a monochromatic configuration, $X_t$ is equally likely to be any of the chambers with the same number of blue vertices by symmetry. Thus, we have
\begin{align*}
\sum_{C \in \mathcal{C}_k} |P(X_t=C) - \pi(C) | & = \left| \sum_{C \in \mathcal{C}_k} P(X_t=C) - \pi(C)   \right| \\
& = | \mathbf{P}(\Phi(X_t)=k) - \pi_{\Phi}(k) |.
\end{align*}

The equality \eqref{egalite} shows that we can study the convergence of $\Phi(X_t)$ to its stationary distribution to obtain an upper bound for the mixing time of the edge flipping on $K_n$. We define a coupling with three stages to study the mixing time. Let $X_t$ and $Y_t$ be two edge flipping process where the initial states are monochromatic of opposite colors. We denote the number of blue vertices in each process by $W_t=\Phi(X_t),Z_t=\Phi(Y_t),$ and take $W_0=0$ and $Z_0=n.$ From now on, we will only deal with $W_t$ and $Z_t$. Let us first define
\begin{equation*}
\Delta_t:=Z_t-W_t \textnormal{ for } t=0,1,\ldots
\end{equation*}

\noindent \textit{First stage: Coupon collector's problem} 

For the first stage of the coupling $(W_t,Z_t)_{t=0}^{\infty},$ we define the stopping time $\tau_1$ as  
\begin{equation*}
\tau_1=\min_t \{ \Delta_t \leq \lceil\sqrt{n}\rceil \},
\end{equation*}
then we define the rule of coupling up to this time. If $t < \tau_1,$ the two chains move identically in the sense that we choose the same two vertices in the underlying graphs and color them the same color. 
 
The first part suggests that the process up to time $\tau_1$ is as twice as fast as the coupon collector's problem where you stop collecting more coupons if the remaining number of coupons is $\lceil\sqrt{n}\rceil.$ Let us define a process $\tilde{\Delta}_t$ for all $t$ with transition probabilities:
 \[\tilde{\Delta}_{t+1} = \tilde{\Delta}_t + \begin{cases} 
      0 & \textnormal{ with prob. } \binom{n-\tilde{\Delta}_t}{2} / \binom{n}{2} \\
      -1 & \textnormal{ with prob. } \tilde{\Delta}_t(n-\tilde{\Delta}_t) \mathbin{/} \binom{n}{2} \\
      -2 & \textnormal{ with prob. } \binom{\tilde{\Delta}_t}{2} \mathbin{/} \binom{n}{2}
   \end{cases}
\]
and $\tilde{\Delta}_0=n.$ Note that $\tilde{\Delta}_t=\Delta_t$ for $t \leq \tau_1.$ Calculating the conditional expectation, we have
\begin{equation*}
\mathbf{E}[\tilde{\Delta}_{t+1}  | \tilde{\Delta}_t] = \left(1-\frac{2}{n} \right)\tilde{\Delta}_t.
\end{equation*}
Therefore,
\begin{equation*}
\mathbf{E}[\tilde{\Delta}_t ] = \left(1-\frac{2}{n}\right)^t \tilde{\Delta}_0 \leq n e^{-\frac{2t}{n}}.
\end{equation*}
By Markov's inequality,
\begin{equation}
\mathbf{P}(\tau_1 > s )=\mathbf{P}[\tilde{\Delta}_s > \sqrt{n} ]  \leq \sqrt{n} e^{-\frac{2s}{n}}.
\end{equation}
If we take $s=\frac{1}{4}n \log n+\gamma_1 n,$ we obtain
\begin{equation}\label{couple1}
\mathbf{P}(\tau_1 > \frac{1}{4}n \log n+\gamma_1 n ) < e^{-\gamma_1}.
\end{equation}

\noindent \textit{Second stage: Coupling with lazy simple random walks} 
 
Note that if we continued with the first coupling until the two chains meet, this expression would give an upper bound of order $\frac{1}{2}n \log n,$ which is the same with what we obtained from eigenvalues. To improve this rate, whenever $t\geq \tau_1,$ we let the chains move according to the second rule where we choose edges independently for two graphs but color them with the same color. We will show that it only needs $\mathcal{O}(n)$ steps to hit zero with positive probability. Given that $t \geq \tau_1,$

\begin{equation} \label{delta}
\Delta_{t+1}=\Delta_t + \begin{cases} 
      -2 & \textnormal{ with prob. } \binom{n-W_t}{2}\binom{Z_t}{2} / \binom{n}{2}^2  \\
      -1 & \textnormal{ with prob. } \left[\binom{n-W_t}{2} Z_t(n-Z_t) + W_t(n-W_t)\binom{Z_t}{2} \right]  / \binom{n}{2}^2  \\  
       1 & \textnormal{ with prob. } \left[W_t(n-W_t)\binom{n-Z_t}{2} + \binom{W_t}{2}Z_t(n-Z_t)\right]  / \binom{n}{2}^2  \\      
       2 & \textnormal{ with prob. } \binom{W_t}{2}\binom{n-Z_t}{2} / \binom{n}{2}^2 \\
             0 & \textnormal{ otherwise. }  \\ 
   \end{cases}
\end{equation}
 
We need to be more precise and show that the terms in the difference $\Delta_t$ are concentrated around $pn$ recalling that $p$ was the probability of coloring a chosen edge blue. 
\begin{lemma}\label{concentration}
 If $t > \frac{1}{4} n \log n,$ then $\mathbf{P}(|V_t-pn| >n^{1/2+\alpha}) \leq 2n^{-2\alpha} (1+o(1))$ for any $\alpha > 0$ where $V_t$ can be replaced by $W_t$ or $Z_t,$ which are defined in the beginning of the proof. 
\end{lemma}
\pf It is not difficult to see that $\mathbf{E}[W_t] \rightarrow pn$ as $t \rightarrow \infty$. To show the concentration result, we will use a second-order estimate. First, consider
\begin{equation}\label{subsupmart}
 W_{t+1}=W_t + \begin{cases} 
      -2 & \textnormal{ with prob. } q  \binom{W_t}{2} / \binom{n}{2}  \\
      -1 & \textnormal{ with prob. } q   W_t(n-W_t) / \binom{n}{2}  \\   
       1 & \textnormal{ with prob. } p W_t(n-W_t) / \binom{n}{2}  \\      
       2 & \textnormal{ with prob. } p  \binom{n-W_t}{2} / \binom{n}{2}\\
             0 & \textnormal{ otherwise }  \\ 
   \end{cases}
\end{equation}
where we take $\binom{z}{n}=0$ unless $z$ is a positive integer. Let 
\begin{equation}\label{sigmaf}
\mathcal{F}_t=\sigma(W_{0}, W_1,\ldots, W_t)
\end{equation}
 be the $\sigma$-field generated by $W_{0},W_{1},\ldots, W_t.$ The computations give $\mathbf{E}[W_{t+1}|\mathcal{F}_t]=2p+\left( 1-\frac{2}{n}\right)W_t.$ We can evaluate this sum recursively to obtain
\begin{equation}\label{recexpect}
\mathbf{E}[W_t]=pn -\left(1-\frac{2}{n}\right)^{t+1}pn > pn - p\sqrt{n}
\end{equation}
provided that $t > \frac{1}4 n \log n$ and $W_0=0$ as assumed earlier. 

Secondly, we bound the variance of $W_t$ by applying a crucial step in the proof of Wilson's Lemma \ref{Wilson} (see Section 13.2 of \cite{L}), which gives

\begin{equation}\label{varbound}
\textnormal{Var}(W_t) \leq \frac{R}{1-\lambda} = \frac{4}{1-\left(1-\frac{2}{n}\right)}=2n
\end{equation}
by \eqref{rbound}. Then by \eqref{recexpect} and using \eqref{varbound} in Chebyshev's inequality,
\begin{equation}
\mathbf{P}\left(|W_t-pn| \geq (p+r\sqrt{2})\sqrt{n} \right) \leq \frac{1}{r^2}
\end{equation}
for $r>0.$ This proves the result for $V_t=W_t.$ The proof for $V_t=Z_t$ follows exactly the same line of argument. By symmetric replacement of parameters, we have 

\begin{equation*}
\mathbf{E}[Z_t]=pn -\left(1-\frac{2}{n}\right)^{t+1}pn + \left(1-\frac{2}{n}\right)^{t} Z_0  < pn +q\sqrt{n}
\end{equation*}
where we take $Z_0=n.$ Since the bound in \eqref{rbound} is true for all initial states, we have the same variance bound. The rest follows from Chebyshev's inequality. \bbox

Next we will show a maximal inequality for $W_t$ in the second stage of the coupling.

\begin{lemma}\label{maximal} 
If $t\geq \frac{1}{4}n \log n$ and $k$ is a positive real number independent of $n$, then\\ $\mathbf{P}\left(\max\limits_{\tau_1 \leq t \leq \tau_1+kn }|V_t -pn| > n^{1/2+\beta}\right) \leq 4n^{-\alpha}$ for any $\beta>0$ and $\alpha \in (0,\beta)$ where $V_t$ can be replaced by $W_t$ or $Z_t.$ 
\end{lemma}

\pf Since $W_t$ and $Z_t$ are symmetric cases as in Lemma \ref{concentration}, we only prove the result for $W_t.$ By Doob-Meyer decomposition, we construct a martingale from $\{W_{\tau_1+s}\}_{s \geq 0}$ with respect to the $\sigma$-fields:
\begin{equation}
\mathcal{G}_s:=\mathcal{F}_{\tau_1+s} \quad \tn{ for all } \quad s=0,1,2,\ldots
\end{equation}
where $\mathcal{F}_t$ is as defined in \eqref{sigmaf} for all $t=0,1,2,\ldots$  Since $\mathbf{P}(\tau_1 < \infty)=1$ by \eqref{couple1}, $M_1=W_{\tau_1+1}-W_{\tau_1}$ is well-defined. For $s\geq 2,$ we let
\begin{equation}\label{martin}
M_s=W_{\tau_1+s} +\frac{2}{n}\sum_{i=1}^{s-1}W_{\tau_1+i} -W_{\tau_1} -2sp.
\end{equation}
From \eqref{subsupmart} and \eqref{recexpect}, we can prove the martingale property:
\[\mathbf{E}[M_{s+1}|\mathcal{G}_{s}]=M_s.\]
Observe that the transition probabilities of the martingale depend on $W_{\tau_1}.$ We will fix $W_{\tau_1}$ later in the proof.

In order to bound the variance of the martingale, we consider the variance of the differences. The computations give
\begin{equation*}
\mathbf{E}[(M_{s+1}-M_s)^2]=C_1(p)\frac{\mathbf{E}[W_{\tau_1+s}^2]}{n^2}+C_2(p)\frac{\mathbf{E}[W_{\tau_1+s}]}{n}+C_3(p) \leq C(p)
\end{equation*}
 where $C(p)$ and the other constants in the expression are independent of $n$ and $s.$ Therefore, by the orthogonality of martingale differences
\[\textbf{E}[M_s^2]=\textnormal{Var}(M_s)\leq C(p)s.\]
Then, by Doob-Kolmogorov maximal inequality, see Theorem 2 in Section 7.8 of \cite{GS20},
\begin{equation}\label{maxrealised}
\mathbf{P}\left(\max\limits_{1 \leq t \leq kn}|M_t|\geq a\right) \leq \frac{\mathbf{E}[M_{kn}]}{a^2} = \frac{C(p)kn}{a^2}
\end{equation}
for $a>0.$

Next, we show that, for $a$ of order $n^{1/2+\beta}$ and large enough $n,$ $|W_{\tau_1+s}-pn|\geq a$ implies $|M_{s}|\geq n^{-\beta/2}a$ (not optimal) with high probability. Let $A$ be the event that 
\begin{equation}\label{lem21}
|W_{\tau_1}-pn|\leq n^{1/2+\beta/2},
\end{equation}
which has probability greater than $1-4n^{-\beta}$ by Lemma \ref{concentration}. The following argument is conditioned on $A.$ Let us fix an outcome $W_{\tau_1}=W^*$ in $A.$ Note that both the stopping time  $\tau_1$ and $W_{\tau_1}$ are determined. We take $n > 2^{4/\beta}$. Suppose 
\begin{equation}\label{lem22}
|W_{\tau_1+s}-pn| > n^{1/2+\beta}.
\end{equation}
From \eqref{lem21} and \eqref{lem22}, we get 
\[|W_{\tau_1+s}-W_{\tau_1}| > \frac{1}{2}n^{1/2+\beta}.\] 
There are two cases for \eqref{lem22}, namely $W_{\tau_1+s} > pn+n^{1/2+\beta}$ or $W_{\tau_1+s} < pn-n^{1/2+\beta}.$ Let us assume the former. It will be apparent below that the other case can be treated similarly. From \eqref{martin}, we have either 
\begin{equation}\label{lem24}
M_s > \frac{1}{4}n^{1/2 + \beta}
\end{equation}
or 
\begin{equation}\label{lem25}
\frac{2}{n}\sum_{i=1}^{s-1}W_{\tau_1+i} -2sp \leq -\frac{1}{4}n^{1/2 + \beta}.
\end{equation}
If the latter is the case, by \eqref{lem21} and \eqref{lem22}, there must exist a time $s'$ such that

\[s'= \max \left\{ 0 \leq t < s \, : \, W_{\tau_1+t} \leq pn+2n^{1/2+\beta/2} \right\}.  \]
Then observe that
\[ \frac{2}{n}\sum_{i=1}^{s'-1}W_{\tau_1+i}  -2s'p < \frac{2}{n}\sum_{i=1}^{s-1}W_{\tau_1+i}  -2sp \]
since $pn < W_{\tau_1+t}$ for $s' \leq t <s.$ Therefore,
\begin{align*}
M_{s'} &= W_{\tau_1+s'}-W_{\tau_1}+\frac{2}{n}\sum_{i=1}^{s'-1}W_{\tau_1+i}  -2s'p \\
& \leq 3 n^{1/2+\beta/2} -\frac{1}{4}n^{1/2 + \beta} \\
& \leq - n^{1/2+\beta/2}.
\end{align*}
In either case, we have 
$$\max\limits_{1 \leq t \leq s}|M_t|\geq n^{1/2+\beta/2}.$$ 
Therefore, by \eqref{maxrealised}, 
\begin{align*}
\mathbf{P}\left(\max\limits_{\tau_1 \leq t \leq \tau_1+s }|W_t -pn| > n^{1/2+\beta} \, \big| \, W_{\tau_1}=W^* \right) & < \mathbf{P}\left(\max\limits_{1 \leq t \leq s}|M_t|\geq n^{1/2+\beta/2} \, \big| \, W_{\tau_1}=W^* \right)\\
 &<C(p)kn^{-\beta}
\end{align*}
Now let us take $s=kn.$ Since the established upper bound is uniform for   $W^*$ in $A$, it follows via a conditioning argument  that 

\begin{align*}
\mathbf{P}\left(\max\limits_{\tau_1 \leq t \leq \tau_1+kn }|W_t -pn| > n^{1/2+\beta}\right) & < \mathbf{P}\left(\max\limits_{1 \leq t \leq kn}|M_t|\geq n^{1/2+\beta/2}  \, \big| \, A\right)\mathbf{P}(A) + \mathbf{P}(A^c) \\
 &<C(p)kn^{-\beta}\mathbf{P}(A)+ \mathbf{P}(A^c) \\
& \leq C(p)kn^{-\beta}(1-4n^{-\beta})+4n^{-\beta} \\
& < 4n^{-\alpha}
\end{align*}
where $n$ is chosen large enough to satisfy the last inequality. \bbox

We return to \eqref{delta}. We use the lemmas above to bound the transition probabilities of $\Delta_t$ uniformly. We first need the following corollary. 
\begin{corollary}\label{corollary}
Let $t \in \left[\tau_1, \tau_1+kn \right]$ for some $k>0$ independent of $n.$  For every $\epsilon > 0,$ there exist $\alpha \in (0,1/2)$ and $\beta \in (\alpha,1/2)$ such that
\begin{equation}\label{pbound}
p_t:=\frac{\sqrt{2} W_t}{n} \in [p-\epsilon, p+\epsilon]
\end{equation}
and 
  \[\delta:=\sqrt{2}\frac{Z_t-W_t}{n}=\mathcal{O}(n^{\beta-1/2})\]
with probabilty $1-8n^{-\alpha}$.
\end{corollary}
\pf Both statements follow from Lemma \ref{maximal}. For the latter, since $W_t$ and $Z_t$ are correlated by coupling, we use the union bound. \bbox

Therefore, under the assumptions of Corollary \ref{corollary}, we can rewrite \eqref{delta} as
 \[\Delta_{t+1} = \Delta_t + \begin{cases} 
      -2 & \textnormal{ with prob. } p_t^2q_t^2 +2p_tq_t^2\delta + o(n^{2\beta-1})  \\
      -1 & \textnormal{ with prob. } 2p_tq_t(p_t^2+q_t^2)+q_t((q_t-p_t)^2+p_t)\delta +o(n^{2\beta-1})  \\
             1 & \textnormal{ with prob. } 2p_tq_t(p_t^2+q_t^2)-p_t((p_t-q_t)^2+q_t)\delta         +o(n^{2\beta-1})  \\      
       2 & \textnormal{ with prob. }  p_t^2q_t^2 -2p_t^2q_t\delta + o(n^{2\beta-1}) \\
              0 & \textnormal{ otherwise. }  \\ 
   \end{cases}
\]
 Then, examining the probabilities, we see that it is more likely for $\{\Delta_{t}\}_{t > \tau_1}$ to go to left if $\Delta_t >0.$ In fact, 
\begin{align}\label{comparison}
\begin{split}
\mathbf{P}(\Delta_{t+1}-\Delta_t =-2)& \geq p_t^2q_t^2\geq \mathbf{P}(\Delta_{t+1}-\Delta_t =2)  \\\
 \mathbf{P}(\Delta_{t+1}-\Delta_t =-1) &\geq 2p_tq_t(p_t^2+q_t^2) \geq \mathbf{P}(\Delta_{t+1}-\Delta_t =1) 
 \end{split}
\end{align}
provided that $Z_t\geq W_t.$ Thus, we can compare it to the following random walk, which is equally likely to go to either direction at any time.
 \[S_{t+1} = S_t + \begin{cases} 
      -2 & \textnormal{ with prob. } p_t^2q_t^2  \\
      -1 & \textnormal{ with prob. } 2p_tq_t(p_t^2+q_t^2)  \\
       0 & \textnormal{ with prob. } (p_t^2+q_t^2)^2+2p_t^2q_t^2 \\ 
       1 & \textnormal{ with prob. } 2p_tq_t(p_t^2+q_t^2)  \\      
       2 & \textnormal{ with prob. }  p_t^2q_t^2 
   \end{cases}
\]
Let us take $S_0=\Delta_{\tau_1}=\lceil \sqrt{n} \rceil.$ We define the following stopping times
\begin{equation*}
\tau_2=\min_t \{\Delta_{\tau_1+t} \in I\},
\end{equation*}
\begin{equation*}
\tau_S=\min_t \{S_t \in I\}
\end{equation*}
where $I=\{-1,0,1\}.$ Comparing the transition probabilities of the two random walks by \eqref{comparison}, we have $\Delta_{\tau_1+t}\leq S_t$ since $Z_t \geq W_t.$ Therefore, we can find a coupling for two processes such that 
\begin{equation}\label{tau2taus}
\tau_2 \leq \tau_S. 
\end{equation}
Furthermore, observe that $S_{t+1}-S_t$ is sum of the i.i.d. random variables which take values $-1$ or $1$ with probability $p_tq_t$ and $0$ with probability $p_t^2+q_t^2.$ But then $S_t=S_{2t}'$ where $S'_{t}$ is none other than a lazy simple random walk on integers, so we may use the results on the lazy simple random walk to bound $\mathbf{P}(\tau_S > s)$, such as Corollary 2.28 of \cite{L}. The only difficulty is that $S_t$ is not time-homogenous as the parameter $p_t$ varies, which is to be dealt with next. 

Given $r \in [0,1],$ we let
 \[\tau_S(r)=\min_t\{S_t(r) \in I\}\]
 where $S_t(r)$ is the lazy simple random walk with holding probability $r^2+(1-r)^2.$ Considering the bound \eqref{pbound}, we define 
 \[p^{*}= (p-\epsilon)\mathbf{1}_{\{p\leq 1/2\}}+ (p+\epsilon)\mathbf{1}_{\{p > 1/2\}}\] 
 The choice of $p^*$ ensures that $S_t(p^*)$ has the maximal holding probability among all values in $[p-\epsilon,p+\epsilon].$  


\begin{lemma}\label{pminuseps}
$\mathbf{P}(\tau_S(p^*) > s) \geq \mathbf{P}(\tau_S > s).$
\end{lemma}
\pf Let us take $q^*=1-p^*.$ We define a Bernoulli random variable for each $t$ as follows:
\[\xi_t = \begin{cases}
0 & \textnormal{ with prob. } \frac{1-(p^{*2}+q^{*2})}{1-(p_t^2+q_t^2)}  \\
      1 & \textnormal{ otherwise. }  \\
\end{cases}\]

$\xi$ is well-defined by the choice of $p^*.$ Now observe that $S_t(p^{*})$ has the same distribution with the chain which moves according to the law of $S_t$ (which is the same with $S_t(p_t)$) if $\xi_t=1,$ stays at the same state if $\xi_t=0.$ Then, it is obvious that we can find a coupling satisfying $\tau_S(p^*) \geq \tau_S.$ The result follows. \bbox

Finally, by Lemma \ref{maximal}, \ref{pminuseps}, Eq. \eqref{tau2taus} and modifying the numbers in Corollary 2.28 of \cite{L} accordingly, we have
\begin{equation*}
\mathbf{P}(\tau_2 > s) \leq 4n^{-2\alpha}+\frac{4\sqrt{n}}{\sqrt{2p^*(1-p^*)s}}
\end{equation*}
So, if we take $s=\gamma_2 n,$ 
\begin{equation}\label{couple2}
\mathbf{P}(\tau_2 > \gamma_2 n) \leq 4n^{-2\alpha}+ \frac{4}{\sqrt{2p^*(1-p^*)\gamma_2}} \leq \frac{4 \max\{p^{-1/2}, (1-p)^{-1/2}\}}{\sqrt{\gamma_2}}
\end{equation}
by choosing $\epsilon$ small enough, let us say $\epsilon=\frac{\min\{p,1-p\}}{100}$ in \eqref{pbound}.\\

\noindent \textit{Third stage: Number of sign changes}

For the last stage of the coupling, let us define 
\begin{equation*}
\tau_3=\min_t \{ \Delta_{\tau_1 +\tau_2+t} = 0 \},
\end{equation*}
which is to account for the case that $Z_{\tau_2}-W_{\tau_2}=1$ and $\Delta_t$ does not hit zero in the next step. By symmetry, we see that it is more likely for $\{\Delta_{t}\}_{t > \tau_1}$ to go to right if $\Delta_t <0.$ So, we again can bound the stopping time by the time associated with $S_t.$ In fact, for the simple random walk, we have the probability that ``the number of sign changes is equal to k up to $t=2n+1$" is equal to ``the probability that the walk is at position $2k+1$ at $t=2n+1,$" see Theorem 1 in Chapter III.6 of \cite{F68}. This allows us to use the normal approximation for the sign changes (Theorem 2  in Chapter III.6 of \cite{F68}). Let $F_Z$ be the distribution function of the standard normal distribution. Let us use $S_t(p^*)$ again, the laziest simple random walk in the range. We have 
\begin{equation*}
\mathbf{P} \left(|t \in \{1,2,\ldots,K(p^*)n \} \, : \, S_{\tau_1+\tau_2+t}(p^*) \in I| \leq \frac{\sqrt{n}}{\gamma_3}\right) \leq 2F_Z\left(\frac{1}{\gamma_3}\right)-1 \leq \frac{1}{\sqrt{2\pi} \gamma_3}.
\end{equation*}
where $K(p)=1+\frac{p^2+q^2}{2pq}.$ Since
\begin{equation*}
\mathbf{P} \left(|t \in \{1,2,\ldots,k \} \, : \, S_{\tau_1+\tau_2+t} \in I| \leq x\right)
 \leq\mathbf{P} \left(|t \in \{1,2,\ldots,k \} \, : \, S_{\tau_1+\tau_2+t}(p^*) \in I| \leq x\right)
\end{equation*}
for all $k$ and $x$, $\Delta_{t}$ visits $I$ of order larger than $\sqrt{n}$ times a probability bounded away from zero. At every visit to $I,$ it has a positive probability to hit zero at the next stage. Let $\kappa=\mathbf{P}(S_{t+1}(p^*)=0 \,| \, S_t=-1 \text{ or }1).$ Thus, 
\begin{equation}\label{couple3}
\mathbf{P}(\tau_3 > K(p^{*}) n) \leq (1-\kappa)^{\frac{\sqrt{n}}{\gamma_3}}+ \frac{1}{\sqrt{2\pi} \gamma_3} < \gamma_3^{-1}
\end{equation}
for any $\gamma_3>0$ independent of $n$ for large $n.$ Let us take $\gamma_3>K(p^*)$ for convenience. Finally let $\tau=\tau_1+\tau_2+\tau_3,$ which is the stopping time, $\alpha>0$ and $\gamma=\gamma_1+\gamma_2+\gamma_3.$ By \eqref{couple1}, \eqref{couple2},  \eqref{couple3},   and using the strong Markov property iteratively,
\begin{equation*}
\mathbf{P}\left(\tau> \frac{1}{4}n \log n + \gamma n\right) <1-\left(1-e^{-\gamma_1}\right)\left(1-4 \max\{p^{-1/2}, (1-p)^{-1/2}\}\gamma_2^{-1/2}\right)\left(1-\gamma_3^{-1}\right).
\end{equation*}

Since we chose the opposite monochromatic configurations for the different initial states of the coupled chains, among all choices for the starting pair chambers the coupling time above is clearly the maximum, so is the probability in \eqref{coupletime}. Therefore,
\begin{equation}
d\left(  \frac{1}{4}n \log n + \gamma n\right) \leq \frac{C}{\sqrt{\gamma}} 
\end{equation}
by choosing $\gamma_1, \gamma_2$ and $\gamma_3$ large enough. Combining with the lower bound in Theorem \ref{epthm}, we have $t_n=\frac{1}{4} n \log n$ and $w_n=n$ in Definition \eqref{cutoffdef}. \bbox

 \bigskip
 
\noindent \textbf{Acknowledgments:} 
The second author is supported partially by BAP grant 20B06P1. The second and the third authors would like to thank Jason Fulman for the suggestion of this subject of study. We would also like to thank the anonymous referees whose suggestions and corrections improved the paper significantly.

\bibliographystyle{alpha}
\bibliography{references}

\begin{thebibliography}{BCCG15}

\bibitem[AD10]{AD10}
C.~A. Athanasiadis and P.~Diaconis.
\newblock Functions of random walks on hyperplane arrangements.
\newblock {\em Advances in Applied Mathematics}, 45(3):410--437, 2010.

\bibitem[BCCG15]{B15}
S.~Butler, F.~Chung, J.~Cummings, and R.~Graham.
\newblock Edge flipping in the complete graph.
\newblock {\em Advances in Applied Mathematics}, 69:46--64, 2015.

\bibitem[BD98]{BD98}
K.~Brown and P.~Diaconis.
\newblock Random walks and hyperplane arrangements.
\newblock {\em The Annals of Probability}, 26(4):1813--1854, 1998.

\bibitem[BHP17]{BHP17}
R.~Basu, J.~Hermon, and Y.~Peres.
\newblock Characterization of cutoff for reversible {Markov} chains.
\newblock {\em Annals of Probability}, 45(3):1448--1487, 2017.

\bibitem[BHR99]{BHR99}
P.~Bidigare, P.~Hanlon, and D.~Rockmore.
\newblock A combinatorial description of the spectrum for the tsetlin library
  and its generalization to hyperplane arrangements.
\newblock {\em Duke Mathematical Journal}, 99(1):135--174, 1999.

\bibitem[Bro00]{B00}
K.~Brown.
\newblock Semigroups, rings, and {M}arkov chains.
\newblock {\em Journal of Theoretical Probability}, 13(3):871--928, 2000.

\bibitem[CG12]{CG12}
F.~Chung and R.~Graham.
\newblock Edge flipping in graphs.
\newblock {\em Advances in Applied Mathematics}, 48:37--63, 2012.

\bibitem[CH14]{C14}
F.~Chung and J.~Hughes.
\newblock A note on an alternating upper bound for random walks on semigroups.
\newblock {\em Discrete Applied Mathematics}, 176:24--29, 2014.

\bibitem[Fel68]{F68}
W.~Feller.
\newblock {\em Introduction to Probability and its Applications}, volume~17.
\newblock Wiley, 1968.

\bibitem[GS20]{GS20}
G.~Grimmett and D.~Stirzaker.
\newblock {\em Probability and random processes}.
\newblock Oxford University Press, 2020.

\bibitem[LPW09]{L}
D.~A. Levin, Y.~Peres, and E.~L. Wilmer.
\newblock {\em {M}arkov Chains and Mixing Time}.
\newblock AMS, Providence, RI, 2009.

\bibitem[Nes17]{N17}
E.~Nestoridi.
\newblock A non-local random walk on the hypercube.
\newblock {\em Advances in Applied Probability}, 49(4):1288--1299, 2017.

\bibitem[Nes19]{N19}
E.~Nestoridi.
\newblock Optimal strong stationary times for random walks on the chambers of a
  hyperplane arrangement.
\newblock {\em Probability Theory and Related Fields}, 174(3):929--943, 2019.

\bibitem[Pik13]{JP}
J.~Pike.
\newblock {\em Eigenfunctions for random walks on hyperplane arrangements}.
\newblock PhD thesis, University of Southern California, 2013.

\bibitem[Sal12]{Sal}
F.~Saliola.
\newblock Eigenvectors for a random walk on a left-regular band.
\newblock {\em Advances in Applied Mathematics}, 48(2):306--311, 2012.

\bibitem[Sal21]{S21}
J.~Salez.
\newblock Cutoff for non-negatively curved {Markov chains}.
\newblock {\em arXiv preprint arXiv:2102.05597}, 2021.

\bibitem[Sta07]{St07}
R.~P. Stanley.
\newblock An introduction to hyperplane arrangements.
\newblock {\em Geometric Combinatorics}, 13:389--496, 2007.

\bibitem[Wil04]{Wilson}
D.~B. Wilson.
\newblock Mixing times of lozenge tiling and card shuffling {M}arkov chains.
\newblock {\em The Annals of Applied Probability}, 14(1):274--325, 2004.

\end{thebibliography}
 
 \end{document}